\documentclass[12pt]{article}

\usepackage{amssymb}
\font\smallit=cmti10
\font\smalltt=cmtt10

\def\B{\hfill $\Box$}
\def\ni{\noindent}

\begin{document}
\begin{center}
{\uppercase{\bf On Monochromatic Ascending Waves}}

\vskip 20pt
{\bf Tim LeSaulnier\footnote{
This work was done as part of a high honor thesis in mathematics
while the first author was an undergraduate at Colgate University,
under the directorship of the second author.
}}\\
\vskip 5pt
and
\vskip -5pt
{\bf Aaron Robertson}\\
{\smallit Department of Mathematics,
Colgate University,
Hamilton, NY 13346}\\
{\smalltt aaron@math.colgate.edu}
\end{center}

\begin{abstract}
A sequence of positive integers $w_1,w_2,\dots,w_n$ is called an
ascending wave if $w_{i+1}-w_i \geq w_i - w_{i-1}$ for $2 \leq i
\leq n-1$.  For integers $k,r\geq1$, let $AW(k;r)$ be the least
positive integer such that under any $r$-coloring of $[1,AW(k;r)]$
there exists a $k$-term monochromatic ascending wave. The existence
of $AW(k;r)$ is guaranteed by van der Waerden's theorem on
arithmetic progressions since an arithmetic progression is, itself,
an ascending wave. Originally, Brown, Erd\H{o}s, and Freedman
defined such sequences and proved that $k^2-k+1\leq AW(k;2) \leq
\frac{1}{3}(k^3-4k+9)$. Alon and Spencer then showed that $AW(k;2) =
O(k^3)$.  In this article, we show that
$AW(k;3) = O(k^5)$ as well as  offer a proof of the existence of
$AW(k;r)$ independent of van der Waerden's theorem. Furthermore, we
prove that for any $\epsilon > 0$,
$$
\frac{k^{2r-1-\epsilon}}{2^{r-1}(40r)^{r^2-1}}(1+o(1))
\leq AW(k;r)
\leq
\frac{k^{2r-1}}{(2r-1)!}(1+o(1))
$$
holds for all $r \geq 1$,
which, in particular, improves upon the best known upper bound for
$AW(k;2)$. Additionally, we show that for fixed $k \geq 3$,
$$
AW(k;r)\leq\frac{2^{k-2}}{(k-1)!} \,r^{k-1}(1+o(1)).
$$
\end{abstract}

\section*{\normalsize 0. Introduction}

A sequence of positive integers $w_1,w_2,\dots,w_n$ is called an
{\it ascending wave} if $w_{i+1}-w_i \geq w_i - w_{i-1}$ for $2 \leq
i \leq n-1$. For $k,r \in \mathbb{Z}^+$, let $AW(k;r)$ be the least
positive integer such that under any $r$-coloring of $[1,AW(k;r)]$
there exists a $k$-term monochromatic ascending wave. Although
guaranteed by van der Waerden's theorem,  the existence of $AW(k;r)$
can be guaranteed independently, as we will show.

Bounds on $AW(k;2)$ have appeared in the literature.  Brown,
Erd\H{o}s, and Freedman [2] showed that for all $k \geq 1$,
$$\begin{array}{ll}k^2-k+1\leq AW(k;2)\leq
\frac{k^3}{3}-\frac{4k}{3}+3.\end{array}$$  Soon after, Alon and
Spencer [1] showed that for sufficiently large $k$,
$$\begin{array}{ll}AW(k;2)>\frac{k^3}{10^{21}}-\frac{k^2}{10^{20}}-
\frac{k}{10}+4.\end{array}$$

Recently, Landman and Robertson [4] proposed the refinement of the
bounds on $AW(k;2)$ and the study of $AW(k;r)$ for $r\geq3$.  Here,
we offer bounds on $AW(k;r)$ for all $r \geq 1$, improving upon the
previous upper bound for $AW(k;2)$.

\section*{\normalsize 1. An Upper Bound}

To show that $AW(k;r) \leq O(k^{2r-1})$ is straightforward. We will
first show that $AW(k;r) \leq k^{2r-1}$ by induction on $r$ (and
hence prove the existence of $AW(k;r)$ without
appealing to van der Waerden's theorem.)  The case $r=1$ is trivial,
thus, for $r \geq 2$, assume $AW(k;r-1)\leq k^{2r-3}$ and consider
any $r$-coloring of $[1,k^{2r-1}]$.  Set $w_1=1$ and let the color
of $1$ be red.  In order to avoid a $k$-term monochromatic ascending
wave there must exist an integer $w_2\in[2,k^{2r-3}+1]$ that is
colored red, lest the inductive hypothesis guarantee a $k$-term
monochromatic ascending wave of some color other than red
(and we are done).
Similarly, there must be an integer $w_3\in[w_2 + (w_2-w_1),w_2 +
(w_2-w_1)+k^{2r-3}-1]$ that is colored red to avoid a monochromatic
$k$-term ascending wave. Iterating this argument defines a
monochromatic (red) $k$-term ascending wave $w_1,w_2,\dots,w_k$,
provided that $w_k \leq k^{2r-1}$.  Since for $i\geq2$, $w_{i+1}
\leq w_i + (w_i-w_{i-1}) + k^{2r-3}$ we see that
$w_{i+1}-w_i \leq i k^{2r-3}$ for $i \geq 1$.
Hence, $w_k-w_1 = \sum_{i=1}^{k-1} (w_{i+1}-w_i)
\leq \sum_{i=1}^{k-1} ik^{2r-3} \leq k^{2r-1}-1$
and we are done.

In this section we provide a better upper bound.
Our main theorem in this section follows.

\noindent {\bf Theorem 1} For fixed $r \geq 1$,
$$\begin{array}{ll}AW(k;r) \leq
\frac{k^{2r-1}}{(2r-1)!}(1+o(1)).\end{array}$$

We will prove Theorem 1 via a series of lemmas, but first we
introduce some pertinent notation.

\ni
{\bf Notation}
For $k\geq 2$ and $M\geq AW(k;r)$, let  $\Psi^{M}(k;r)$ be
the collection of all $r$-colorings of $[1,M]$.  For
$\psi\in\Psi^{M}(k;r)$, let $\chi_k(\psi)$ be the set
of all monochromatic $k$-term ascending waves under
$\psi$.  For each monochromatic $k$-term ascending wave
$w=\{w_1,w_2,\ldots,w_{k}\}\in \chi_k(\psi)$, define the $i^{\mathrm{th}}$
difference, $d_{i}(w)=w_{i+1}-w_i$, for $1\leq i\leq k-1$. For
$\psi\in\Psi^M(k;r)$, define
$$\delta_k(\psi)=\min\{d_{k-1}(w)|w\in\chi_k(\psi)\},$$i.e., the
minimum last difference over all monochromatic $k$-term ascending
waves under $\psi$.  Lastly, define
$$\Delta^M(k;r)=\max\{\delta_k(\psi)|\psi\in\Psi^M(k;r)\}.$$
These concepts will provide us with the necessary tools to prove
Theorem 1.  \par We begin with a recursive bound on $AW(k;r)$, which
also proves the existence of $AW(k;r)$ without appealing to van der
Waerden's theorem on arithmetic progressions.

\noindent {\bf Lemma 1.1}  For  $k,r \geq 1$, let $M(k;1)=k$,
$M(1;r)=1$, $M(2;r)=r+1$, and define, for $k \geq 3$ and $r \geq 2$,
$$
M(k;r) =  M(k-1;r)+\Delta^{M(k-1;r)}(k-1;r)+ M(k;r-1)-1.  $$ Then,
for all $k,r \geq 1$, $ AW(k;r)\leq  M(k;r).$

\noindent {\it Proof.}  Noting that the cases
$k+r=2,3,$ and $4$ are, by definition, true, we proceed by
induction on
$k+r$ using $k+r=5$ as our basis.
We have $M(3;2)=7$.  An easy calculation shows
that $AW(3;2)=7$.  So, for some $n\geq5$, we assume
Lemma 1.1 holds for all $k,r\geq1$ such that $k+r=n$.  Now, consider
$k+r=n+1$. The result is trivial when $k=1\;\mathrm{or}\;2$, or if $r=1$,
thus we may assume $k\geq 3$ and $r\geq2$.  Let $\psi$ be an $r$-coloring
of $[1,M(k;r)]$.  We will show that $\psi$ admits a monochromatic
$k$-term ascending wave, thereby proving Lemma 1.1.

By the inductive hypothesis, under $\psi$ there must be a
monochromatic $(k-1)$-term ascending wave
$w=\{w_1,w_2,\ldots,w_{k-1}\}\subseteq[1, M(k-1;r)]$
with  $d_{k-2}(w)\leq\Delta^{M(k-1;r)}(k-1;r) $.  Let
$$N=[w_{k-1}+\Delta^{M(k-1;r)}(k-1;r),
w_{k-1}+\Delta^{M(k-1;r)}(k-1;r)+ M(k;r-1)-1].$$  If there exists
$q\in N$ colored identically to $w$, then $w\cup\{q\}$ is a
monochromatic $k$-term ascending wave, since
$q-w_{k-1}\geq\Delta^{M(k-1;r)}(k-1;r)\geq d_{k-2}(w)$.  If there is
no such $q\in N$, then $N$ contains integers of at most $r-1$
colors. Thus, since $|N|= M(k;r-1)$, the inductive hypothesis
guarantees that we have a monochromatic $k$-term ascending wave in
$N$.  As
$$w_{k-1}+\Delta^{M(k-1;r)}(k-1;r)+ M(k;r-1)-1\leq M(k;r),$$ this
completes the proof.  \hfill$\Box$

\noindent {\bf Corollary 1.2} Let $k\geq 3$ and $r\geq 2$.  Let
$M(k;r)$ be as in Lemma 1.1.  Then
$$\Delta^{M(k;r)}(k;r)\leq\Delta^{M(k-1;r)}(k-1;r)+ M(k;r-1)-1.$$

\noindent
{\it Proof.}
Let $\psi$, $w$, and $N$ be as defined in the proof of Lemma 1.1.
 If there exists $q\in N$ colored identically to $w$,
then
$$\delta_k(\psi) \leq d_{k-1}(w\cup\{q\})\leq\Delta^{M(k-1;r)}(k-1;r)+
M(k;r-1)-1.$$ If there is no such $q\in N$, then there exists a
monochromatic $k$-term ascending wave, say $v$, in $N$.  Hence, $
\delta_k(\psi) \leq d_{k-1}(v)\leq M(k;r-1)-(k-1)$. Since  $\psi$
was chosen arbitrarily, it follows that
$$\Delta^{M(k;r)}(k;r) \leq \Delta^{M(k-1;r)}(k-1;r)+ M(k;r-1)-1.$$
\hfill$\Box$

 Now, as it is easily seen that $\Delta^{M(2;r)}(2;r)=r$ for all $r\geq1$,
Lemma 1.1 and Corollary 1.2 can be iterated to provide bounds on
$M(k;r)$ for any $k\geq3$ and $r\geq2$.  The following lemma will
provide a means for recursively bounding $M(k;r)$, and thus
$AW(k;r)$, by a function of $k$ for any $r\geq2$.

\noindent {\bf Lemma 1.3} Let $k\geq 3$ and $r \geq 2$.  Let
$M(k;r)$ be as in Lemma 1.1. Then
\renewcommand{\arraystretch}{1.5}
$$
\begin{array}{ll}
M(k;r)&\leq \sum_{i=0}^{k-3}((i+1)
M(k-i;r-1))-\frac{k^2}{2}+\frac{3k}{2}+(k-1)r
\end{array}
$$

\noindent {\it Proof.}  We proceed by induction on $k$.
Consider $ M(3;r)$.  We have
$$ M(3;r)=M(2;r)+\Delta^{M(2;r)}(2;r)+ M(3;r-1)-1.$$
Since $M(2;r)=r+1$ and $\Delta^{M(2;r)}(2;r)=r$,  we have
$$\begin{array}{ll}
M(3;r)= M(3;r-1) +2r=
M(3;r-1)-\frac{3^2}{2}+\frac{3(3)}{2}+2r,\\
\end{array}$$ thereby finishing the case $k=3$ and arbitrary $r$.  Now
assume that Lemma 1.3 holds for some $k\geq 3$.  Lemma 1.1, the
inductive hypothesis, and Corollary 1.2 give us
\renewcommand{\arraystretch}{1.5}
$$
\begin{array}{ll}
M(k+1;r)&= M(k;r)+\Delta^{M(k;r)}(k;r)+
M(k+1;r-1)-1\\
&\leq\sum_{i=0}^{k-3}((i+1)
M(k-i;r-1))-\frac{k^2}{2}+\frac{3k}{2}+(k-1)r\\
&\hskip 40pt+\Delta^{M(k;r)}(k;r)+
M(k+1;r-1)-1\\
&\leq\sum_{i=0}^{k-3}((i+1)
M(k-i;r-1))-\frac{k^2}{2}+\frac{3k}{2}+(k-1)r\\
&\hskip 40pt+\Delta^{M(2;r)}(2;r)$$$$+\sum_{i=0}^{k-3}
 M(k-i;r-1)\\
&\hskip 40pt+ M(k+1;r-1)-(k-2)-1\\
&\leq \sum_{i=0}^{k-2}((i+1)
M(k+1-i;r-1))-\frac{(k+1)^2}{2}+\frac{3(k+1)}{2}+kr\\
\end{array}$$
as desired.  \hfill$\Box$

Now, for $r\geq2$, an upper bound on $M(k;r)$, and thus on
$AW(k;r)$, can be obtained by iterating Lemma 1.3.  We offer one
additional lemma, from which Theorem 1 will follow
by application of Lemma 1.1.

\noindent
{\bf Lemma 1.4}  For $k\geq3$ and $r\geq1$, there exists a
polynomial $p_{r}(k)$ of degree at most $2r-2$ such that
$$\begin{array}{ll}M(k;r) \leq
\frac{k^{2r-1}}{(2r-1)!}+p_r(k).\end{array}$$

\noindent {\it Proof.}  We have $ M(k;1)=k$, so we can take
$p_1(k)=1$, having degree $0$.  We proceed by induction on $r$.
Assume Lemma 1.4 holds for some $r\geq1$ so that
$M(k;r)\leq\frac{k^{2r-1}}{(2r-1)!}+p_r(k)$, where $p_r(k)$ is a
polynomial of degree at most $2r-2$.  Lemma 1.3 gives
$$
\begin{array}{ll}
M(k;r+1)&\leq\sum_{j=3}^{k} ((k-j+1)
M(j;r))-\frac{k^2}{2}+\frac{3k}{2}+(k-1)(r+1)\\
&\leq k\sum_{j=3}^{k}
\left(\frac{j^{2r-1}}{(2r-1)!}+p_r(j)\right)-\sum_{j=3}^{k}
\left((j-1)\left(\frac{j^{2r-1}}{(2r-1)!}+p_r(j)\right)\right)\\
&\hskip 40pt -\frac{k^2}{2}+\frac{3k}{2}+(k-1)(r+1).
\end{array}
$$
By Faulhaber's formula [3], for some polynomial $p_{r+1}(k)$ of
degree at most $2r$, we now have
$$M(k;r+1)\leq
k\frac{\frac{k^{2r}}{2r}}{(2r-1)!}-
\frac{\frac{k^{2r+1}}{2r+1}}{(2r-1)!}+p_{r+1}(k)
=\frac{k^{2r+1}}{(2r+1)!}+p_{r+1}(k)$$
and the proof is complete. \hfill $\Box$

As $AW(k;r)\leq M(k;r)$, Theorem 1 now follows, giving that for
fixed $r\geq1$,$$\begin{array}{ll}AW(k;r) \leq
\frac{k^{2r-1}}{(2r-1)!}(1+o(1)).\end{array}$$ Interestingly, Lemma
1.3 can also be used to show the following corollary.

\noindent {\bf Corollary 1.5}  For fixed $k\geq3$,
$$\begin{array}{ll}AW(k;r)\leq\frac{2^{k-2}}{(k-1)!}\,r^{k-1}(1+o(1)).\end{array}$$

\noindent {\it Proof.}  In analogy to Lemma 1.4, we show that for
$k\geq3$ and $r\geq2$, there exists a polynomial $s_{k}(r)$ of
degree at most $k-2$ such that
$$\begin{array}{ll}M(k;r)\leq
\frac{2^{k-2}}{(k-1)!}\,r^{k-1}+s_{k}(r).\end{array}$$

We proceed by induction on $k$. Let $r\geq2$ be arbitrary.  Lemma
1.3 gives that
$$M(3;r)=M(3;r-1)+2r.$$  As $M(3;1)=3$, we now have a recursive
definition of $M(3;r)$ for all $r\geq2$.  We get, for $r\geq2$,
$$
\begin{array}{ll}M(3;r)=M(3;1)+\sum_{i=2}^{r}2i=r^2+r+1,\end{array}$$ which serves as our
basis. Now, for given $k\geq4$, let
$\hat{s}_{3}(r)=(k-1)r-\frac{k^2}{2}+\frac{3k}{2}$ and assume
Corollary 1.5 holds for all integers $3\leq j\leq k-1$ and for all
$r\geq 2$. Lemma 1.3 yields
$$
\begin{array}{ll}
M(k;r)&\leq\sum_{i=0}^{k-3}((i+1)
M(k-i;r-1))+\hat{s}_3(r)\\
&=M(k;r-1)+\sum_{i=1}^{k-3}((i+1) M(k-i;r-1))+\hat{s}_3(r).
\end{array}
$$
Now, by the inductive hypothesis, for $1 \leq i \leq k-3$, we have
that
$$\begin{array}{ll}M(k-i;r-1)&\leq\frac{2^{k-i-2}}{(k-i-1)!}\,(r-1)^{k-i-1}+s_{k-i}(r-1)\\
&=\frac{2^{k-i-2}}{(k-i-1)!}\,r^{k-i-1}+\tilde{s}_{k-i}(r)\,,\end{array}$$
where $\tilde{s}_{k-i}(r)$ is polynomial of degree at most
$k-i-2\leq k-3$. This gives us that
$$\begin{array}{ll}\sum_{i=1}^{k-3}
(i+1)M(k-i;r-1)+\hat{s}_{3}(r)&\leq\sum_{i=1}^{k-3}
\left((i+1)\left(\frac{2^{k-i-2}}{(k-i-1)!}\,r^{k-i-1}+\tilde{s}_{k-i}(r)\right)\right)+\hat{s}_{3}(r)\\
&=2\cdot\frac{2^{k-3}}{(k-2)!}\,r^{k-2}+\check{s}_{k-1}(r),
\end{array}$$ where $\check{s}_{k-1}(r)$ is a polynomial of degree
at most $k-3$. Hence, we have
$$
\begin{array}{ll}
M(k;r)&\leq M(k;r-1)+2 \cdot \frac{2^{k-3}}{(k-2)!} \,r^{k-2}+\check{s}_{k-1}(r)\\
&=M(k;r-1)+  \frac{2^{k-2}}{(k-2)!} \,r^{k-2}+\check{s}_{k-1}(r).
\end{array}$$
As $M(k;1)=k$, we have a recursive bound on $M(k;r)$ for $r\geq2$.
Faulhaber's formula yields
$$\begin{array}{ll}
M(k;r)\leq
M(k;1)+\sum_{i=2}^{r}\left(\frac{2^{k-2}}{(k-2)!}i^{k-2}+\check{s}_{k-1}(r)\right)
\leq \frac{2^{k-2}}{(k-1)!}r^{k-1}+s_{k}(r),\end{array}
$$
where $s_{k}(r)$ is a polynomial of degree at most
$k-2$.  By Lemma 1.1, the result follows. \hfill$\Box$

\section*{\normalsize 2. A Lower Bound for more than Three Colors}

We now provide a lower bound on $AW(k;r)$ for
arbitrary fixed $r \geq 1$.  We generalize
an argument of Alon and Spencer [1]
 to provide
our lower bound.

We will use $\log(x) = \log_2(x)$ throughout.
Also, by $k=x$ for $x \not \in \mathbb{Z}^+$
we mean $k = \lfloor x \rfloor$.

\noindent
{\bf Theorem 2} For fixed $r \geq 1$ and any
$\epsilon>0$, for $k$ sufficiently large,
$$
AW(k;r) \geq\frac{k^{2r-1-\epsilon}}{2^{r-1} (40r)^{r^2-1}}.
$$

The result is trivial for $r=1$ (since $AW(k;1)=k$).  We will assume
that for $r \geq 2$ the inequality holds for $r-1$ and show
that it holds for $r$.

We proceed by defining a certain type of random
coloring.  To this end,
consider the
$r \times 2r$ matrix  $A_0=(a_{ij})$:
$$
\left[
\begin{array}{ccccccccc}
0&0&1&1&2&2&\dots&(r-1)&(r-1)\\
0&1&1&2&2&3&\dots&(r-1)&0\\
0&2&1&3&2&4&\dots&(r-1)&1\\
\vdots&&&&&\vdots&&\vdots&\vdots\\
0&(r-1)&1&0&2&1&\dots&(r-1)&(r-2)\\
\end{array}
\right]
$$
where, for $j \in \left[0,r-1\right]$,
we have $a_{i,(2j+1)} = j$,  for all $1 \leq i \leq r$,
and
$
a_{i,2j+2} = i+j-1 \, (\bmod r).
$

Next, we define $A_j = A_0 \oplus \bf{j}$ where
$\oplus$ means addition modulo $r$ and
$\bf{j}$ is the $r \times 2r$ matrix with
all entries equal to $j$.

Consider the $r^2 \times 2r$ matrix
$
A =
[A_0\,\,A_1\,\,A_2\,\, \dots\,\,A_{r-1}]^t.
$

In the sequel, we will use the following notation.

\ni
{\bf Notation}  For $r \geq 1$, let
$$
N_r = \frac{1}{2^{r-1}(40r)^{r^2-1}}.
$$

Fix $\epsilon>0$.
Let
$$b = AW\left(\frac{k}{10(4r-4)};r-1\right) - 1,$$
so that by the inductive hypothesis, we have
$$b > N_{r-1}
\left(\frac{k}{10(4r-4)}\right)^{2r-3-\epsilon/2}$$ for $k$ sufficiently
large. Using the colors $0,1,\dots,r-1$, let
$\gamma_i$ be an $(r-1)$-coloring of $b$ consecutive
integers with no monochromatic $\frac{k}{10(4r-4)}$-term
ascending wave, where the $r-1$ colors used are
$\{0,1,\dots,i-1,i+1,i+2,\dots,r-1\}$ (i.e., color
$i$ is not used, and hence the subscript on $\gamma$).

We next describe how we randomly $r$-color  $[1,M_\epsilon]$, where
$$M_\epsilon=N_rk^{2r-1 - \epsilon}.$$

We partition the interval $[1,M_\epsilon]$  into consecutive
intervals of length $b$ and denote the $i^{\mathrm{th}}$ such
interval by $B_i$ and call it a {\it block} (note that  the last
block may be a partial block). Let $C_i$ be the $2r$ consecutive
blocks $B_{2r(i-1)+1}, B_{2r(i-1)+2}, \dots,B_{2ri}$, $i=1,2,\dots,
\lceil\frac{M_\epsilon}{2rb}\rceil$. For each $C_i$, we randomly
choose a row in $A$, say $(s_1,s_2,\dots,s_{2r})$.  We color the
$j^{\mathrm{th}}$ block of $C_i$ by $\gamma_{s_j}$.
 By $col(B_i)$ we mean
the coloring of the $i^{\mathrm{th}}$ block, $1 \leq i \leq
\lceil\frac{M_\epsilon}{b}\rceil$, which is one of
$\gamma_0,\gamma_1,\dots,\gamma_{r-1}$.

The following is immediate by construction.

\noindent
{\bf Lemma 2.1}

(i) For all $1 \leq i \leq 2rb$, $P(col(B_i)=\gamma_c)=\frac{1}{r}$
for each $c = 0 ,1, \dots,r-1$.

(ii) For any $i$, $P(col(B_i) = \gamma_c \,\, \mathrm{ and }\,\,
col(B_{i+1})=\gamma_d) =
\frac{1}{r^2}$ for any $c$ and $d$.

(iii)  The colorings of blocks whose pairwise distances
are at least $2r$ are mutually \vskip -5pt
\hskip 26pt independent.

The approach we take, following
Alon and Spencer [1], is to show that there
exists a coloring such that for any monochromatic
$\frac{k}{2}$-term ascending wave
$w_1,w_2,\dots,w_{k/2}$ we have
$w_{k/2}-w_{k/2-1} \geq O(k^{2r-2-\epsilon/2})$.
The following definition
and  lemma, which are generalizations of
those found in [1], will give us the desired result.

\noindent
{\bf Definition}  An arithmetic progression $x_1<x_2<\cdots<x_t$
is called a {\it good progression} if for each
$c \in \{0,1,\dots,r-1\}$, there exists
 $x_i \in B_j$ such that $col(B_j) = col(B_{j+1})=\gamma_c$.
An arithmetic progression that is not good is
called a {\it bad progression}.

\noindent
{\bf Lemma 2.2} For $k,r \geq 2$, let $t =
\frac{(4r-2)(2r+1)}{\log(r^2/(r^2-1))}
\log (k ) + \frac{(2r+1)(\log(r)+1)}{\log(r^2/(r^2-1))}$.
For $k$ sufficiently large, the probability that there is a bad
progression in a random coloring of $[1,M_\epsilon]$ with difference
greater than $b$ of $t$ terms is at most $\frac12$.

\noindent
{\it Proof.} Let $x_1<x_2<\cdots<x_t$ be a bad progression
with $x_2-x_1>b$.
Then no 2 elements belong to the same block.  For each
$i$, $1 \leq i \leq \frac{t}{2r+1}$, let $C_i$ be
the block in which $x_{(2r+1)i}$ resides, and let
$D_i$ be the consecutive block.  Then, the probability
that the progression is bad is at most
\renewcommand{\arraystretch}{1.5}
$$
p=\sum_{j=1}^r P\left(  \nexists \, i \in \left[1,\frac{t}{2r+1}\right]:
Col(C_i) = Col(D_i)=\gamma_j\right).
$$
We have
$$
\begin{array}{ll} \displaystyle
p&\leq r P\left(  \nexists \, i \in \left[1,\frac{t}{2r+1}\right]:
Col(C_i) = Col(D_i)=\gamma_0\right)\\
&= r \left(\frac{r^2-1}{r^2}\right)^{\frac{t}{2r+1}}\\
&\leq r
\left(\frac{r^2-1}{r^2}\right)^{\frac{(4r-2)}{\log(r^2/(r^2-1))}
\log (k )+ \frac{\log(r)+1}{\log(r^2/(r^2-1))}}\\
&\leq \frac{2^{-1}}{k^{4r-2}}
\end{array}
$$
for $k$ sufficiently large.

Since the number of $t$-term arithmetic progressions
in $[1,M_\epsilon]$ is less than $M_\epsilon^2< k^{4r-2}$,
the  probability that there is a bad progression
is less than
$$
k^{4r-2} \cdot \frac{2^{-1}}{k^{4r-2}} = \frac12,
$$
thereby completing the proof.
\hfill$\Box$

\noindent
{\bf Lemma 2.3}  Let $C$ be an $r$-coloring of
$[1,M_\epsilon]$ having no bad progression with difference greater than
$b$ of $t$ terms ($t$ from Lemma 2.2). For any $\epsilon>0$, for $k$
sufficiently large, any monochromatic $\frac{k}{2}$-term ascending wave
$w_1,w_2,\dots, w_{k/2}$ has
$
w_{k/2}-w_{k/2-1} \geq bk^{1-\epsilon/2}  = O(k^{2r-2-\epsilon/2})
$.

\noindent {\it Proof.}  At most $4r-4$ consecutive blocks can have a
specific color in all of them. (To achieve this say the color is
$0$. The random coloring must have chosen row $1$ followed by row
$r+1$, to have $\gamma_0\gamma_0\gamma_1\gamma_1\cdots\gamma_{r-1}
\gamma_{r-1}\gamma_1\gamma_1\gamma_2\gamma_2\cdots\gamma_0\gamma_0$.)
Since each block has a monochromatic ascending wave of length at
most $\frac{k}{10(4r-4)}-1$, any $4r-4$ consecutive blocks
contribute less than $\frac{k}{10}$ terms to a monochromatic
ascending wave.  After that, the next difference must be more than
$b$.

Let the monochromatic ascending wave be $A=a_1,a_2,\dots,a_{k/2}$.
Then, there exists $i<\frac{k}{10}$ such that $a_{i+1}-a_i\geq b+1$.
Now let $X=x_1,x_2,\dots,x_t$ be a $t$-term good progression with
$x_1=a_{i}$ and $d=x_2-x_1=a_{i+1}-a_i \geq b+1$.

Assume, without loss of generality, that the color of $A$ is $0$.
Since $X$ is a good progression, there exists $x_{j} \in B_\ell$
with $Col(B_\ell)=Col(B_{\ell+1})= \gamma_0$ for some block
$B_\ell$. Since $a_{i+j}\geq x_{j}$ as $A$ is an ascending wave, we
see that $a_{i+j}-a_{i} \geq jd+b+1$. We conclude that $a_{i+t}-a_i
\geq td+b+1$ so that, since $A$ is an ascending wave,
$a_{i+t+1}-a_{i+t} \geq d + \frac{b+1}{t}$. Now, let
$X=x_1,x_2,\dots,x_t$ be the $t$-term good progression with
$x_1=a_{i+t}$ and $d'=x_2-x_1 = a_{i+t+1}-a_{i+t}\geq d
+\frac{b+1}{t} \geq (b+1) \left(1+\frac{1}{t}\right)$.  Repeating
the above argument, we see that $a_{i+2t}- a_{i+t} \geq td'+b+1$ so
that $a_{i+2t}-a_{i+2t-1} \geq d' + \frac{b+1}{t} \geq (b+1)
\left(1+ \frac{2}{t}\right)$. In general,
$$
a_{i+st}-a_{i+st-1} \geq (b+1) \left(1+\frac{s}{t}\right)
$$
for $s=1,2,\dots \frac{2k-5t}{5t}$.  Thus, we have
(with $s= (k^{1-\epsilon/2}-1)t \leq \frac{2k-5t}{5t}$ for
$k$ sufficiently large)
$$
a_{k/2} - a_{k/2-1} \geq (b+1) \left(1+
\frac{(k^{1-\epsilon/2}-1)t}{t}\right) = (b+1)k^{1-\epsilon/2}.
$$
\hfill$\Box$

Using Lemma 2.3, there exists an $r$-coloring of $[1,M_\epsilon]$
such that any $\frac{k}{2}$-term monochromatic ascending wave has
last difference at least  $(b+1)k^{1-\epsilon/2}$. This implies that
the last term of any $k$-term monochromatic ascending wave must be
at least $ \frac{k}{2}+(b+1)k^{1-\epsilon/2} \cdot \frac{k}{2} >
\frac12(b+1)k^{2-\epsilon/2}. $

We have
$$
b+1 \geq N_{r-1} \frac{1}{40^{2r-3 - \epsilon/2}
(r-1)^{2r-3-\epsilon/2}}k^{2r-3-\epsilon/2}
\geq N_{r-1}\frac{k^{2r-3-\epsilon/2}}{40^{2r-1} r^{2r-1}}.
$$
Hence, for $k$ sufficiently large, the last term
of any  monochromatic $k$-term ascending wave
must be greater than
$$
N_{r-1} \cdot \frac{1}{40^{2r-1}
r^{2r-1}}k^{2r-3-\epsilon/2} \cdot \frac{k^{2-\epsilon/2}}{2}
 = N_rk^{2r-1-\epsilon} = M_\epsilon.
$$
Hence, we have an $r$-coloring of $[1,M_\epsilon]$ with no $k$-term
monochromatic ascending wave, for $k$
sufficiently large, thereby proving Theorem 2.

\section*{\normalsize 3. A Lower Bound for Three Colors}

We believe that $AW(k;r) = O(k^{2r-1})$,
however, we have thus far been unable to prove this.
The approach of Alon and Spencer [1], which is to
show that there exists an $r$-coloring
(under  a random coloring scheme) such that every monochromatic
$\frac{3k}{4}$-term ascending wave has
$d_{3k/4-1} > ck^{2r-2}$ does not work for an arbitrary
number of colors with our generalization.  However, for 3 colors, we can
refine their argument to prove that $AW(k;3)=O(k^5)$.

\ni
{\bf Theorem 3}
$$
\frac{k^5}{2^{13} \cdot 10^{39}}\leq AW(k;3) \leq
\frac{k^{5}}{120}(1+o(1))
$$

The upper bound comes from Theorem 1, hence
we need only prove the lower bound.
We use the same coloring scheme as in Section 2
and proceed with a series of lemmas.

Let $w_1,w_2,\dots,w_{k/2}$ be a $\frac{k}{2}$-term
ascending wave.
 From Lemmas 2.2 and 2.3, there exist (many) colorings
such that
for $k$ sufficiently large, $w_{k/2}-w_{k/2-1} > 6b\,\, (=2rb)$.
Before proving Theorem 3, we
introduce the following definition.

\ni
{\bf Definition} We call a sequence $x_1,x_2,\dots,x_n$
with $x_2-x_1 \geq 1$ an
{\it almost ascending wave} if, for $2 \leq i \leq n-1$,
we have $d_i=x_{i+1}-x_i$ with
$d_i \geq d_{i-1}- 1$,
with equality for at least one such $i$ and with
the property that if $d_i=d_{i-1}-1$ and $d_j = d_{j-1}-1$
with $j>i$ there must exist $s$, $i<s<j$, such that
$d_s\geq d_{s-1}+1$.

The upper bound of the following proposition is a slight refinement
of a result of Alon and Spencer [1, Lemma 1.7].

\ni {\bf Proposition 3.1}  Denote by $aw(n)$ the number of ascending
waves of length $n$ with first term given and $d_{n-1} <
\frac{n}{10^{14}}$.   Analogously, let $aaw(n)$ be the number of
almost ascending waves of length $n$ with first term given and
$d_{n-1} \leq \frac{n}{10^{14}}$.  Then, for all $n$ sufficiently
large,
$$
2^{\frac{n}{2}-1}<aw(n)+aaw(n) \leq
2^{\frac{13n}{25}} \cdot \left(\frac{3}{2}\right)^{n/100}.
$$

\ni
{\it Proof.}
We start with the lower bound.
 We start by constructing a sequence
where all of $\frac{n}{2}-1$ slots
contain 2 terms of a sequence.
From a list of $\frac{n}{2}-1$ slots, choose $j$, $0 \leq j \leq
\frac{n}{2}-1$, of them.   In these slots place the
pair $-1,1$.  In the remaining slots put the pair $0,0$.
We now have a  sequence of length $n-2$ or $n-3$.
If the length is $n-2$, put a $2$ at the end;
if the length is $n-3$, put $2,2$ at the end.  We now
have, for each $j$ and each choice of $j$ slots, a distinct
sequence of length
$n-1$.   Denote one such
sequence by
$s_1,s_2,\dots,s_{n-1}$. Now, let $d_1=1$ and $d_i = d_{i-1} + s_{i-1}$
for
$i=2,3,\dots,n$.  Since we have the first term of
an almost ascending, or ascending, wave $w_1,\dots,w_n$ given, such a
wave is determined by its sequence of
differences $w_{i+1}-w_i$.
Above, we have constructed a sequence $\{d_i\}$
of differences that adhere to the rules
of an almost ascending, or ascending, wave.
Hence,
$aw(n;r) + aaw(n;r) > \sum_{j=0}^{\frac{n}{2}-1} {\frac{n}{2}-1 \choose j}
=2^{\frac{n}{2}-1}$.

For the upper bound,
we follow the proof of Alon and Spencer
[1, Lemma 1.7], improving the bound enough to serve our purpose.
Their lemma includes the term
${n + \lceil 10^{-6}n\rceil -1 \choose n-1}$
which we will work on to refine their upper bound
on $aw(n)+aaw(n)$.

First,  we have
$$
{n + \lceil 10^{-6}n\rceil -1 \choose n-1} \leq
{(1+10^{-5})n  \choose n}
$$
for $n$ sufficiently large.

Let $q=(1+10^{-5})^{-1}$, $m=\frac{n}{q}$ and let $H(x)=-x \log(x) -
(1-x)
\log(1-x)$ for $0 \leq x \leq 1$ be the binary entropy function.
Then we have\footnote{Here's a quick derivation:
 For all $n \geq 1$, we have
$\sqrt{2\pi n} e^{1/(12n+1)}(n/e)^n
\leq n!
\leq
\sqrt{2 \pi n} e^{1/(12n)} (n/e)^n$
(see [5]).  Hence,
${m \choose qm} \leq \frac{c}{\sqrt{m(1-q)}} \left(q^{-q}
(1-q)^{-(1-q)}\right)^m$
for some positive $c <e^{-2}$
(so that $\frac{c}{\sqrt{m(1-q)}}<1$
for $m$ sufficiently large).  Using
the base 2 $\log$, this gives
${m \choose qm}
\leq 2^{mH(q)}$.
}
$$
{m \choose qm} \leq 2^{m H(q)}.
$$
Applying this, we have
$$
\begin{array}{ll}
H(q) &= \frac{1}{(1+10^{-5})} \log ((1+10^{-5}))
- \left(\frac{10^{-5}}{(1+10^{-5})}\right)
\log \left(\frac{10^{-5}}{(1+10^{-5})}\right)\\
\end{array}
$$
so that
$$\renewcommand{\arraystretch}{2.0}
\begin{array}{ll}
mH(q) &= \left[ \log (1+10^{-5})
- \frac{1}{10^{5}}
\log \left(\frac{10^{-5}}{1+10^{-5}}\right)\right]n\\
&= \left[\frac{1}{10^{5}}
\log \left(10^{5}(1+10^{-5})^{10^{5}+1}\right)\right]n\\
 &\leq\left[\frac{1}{10^{5}}
\log(e(10^{5}+1))\right]n.
\end{array}
$$

We proceed by noting that
$$
\left[\frac{\log(e(10^{5}+1))}{10^{5}}\right]n
\leq
\left[\frac{1}{100}\log\left(\frac{3}{2}\right)\right]n.
$$
Hence,
$
2^{mH(q)} \leq
2^{\frac{n}{100} \log\left(\frac{3}{2}\right)}
= \left(\frac{3}{2}\right)^\frac{n}{100}.
$
Now, using Alon and Spencer's result [1, Lemma 1.7],
the result follows.
\B

We are now in a position to prove the fundamental
lemma of this section.  In the proof we refer to the following definition.

\ni
{\bf Definition}  Let $a_1,\dots,a_n$ be an ascending
wave and let $b \in \mathbb{Z}^+$.  We call
$\lfloor \frac{a_1}{b} \rfloor, \lfloor \frac{a_2}{b} \rfloor,
\cdots, \lfloor \frac{a_n}{b} \rfloor$ the associated
{\it $b$-floor wave}.

\ni
{\bf Lemma 3.2}  Let
$$
Q=\frac{k^5}{2^{13} \cdot 10^{39}}.
$$
 The probability that in a random $3$-coloring
of $[1,Q]$ there is
a monochromatic ascending wave of length $\frac{k}{4}$ whose
first difference is greater than $6b \,\,(=2rb)$ and whose last difference
is smaller than $\frac{kb}{4 \cdot 10^{14}} = O(k^{4})$
is less than
$\frac{1}{2}$ for $k$ sufficiently large.

\ni
{\it Proof.}
Let $A=a_1<a_2<\cdots<a_{k/4}$ be an ascending
wave and let
$\lfloor \frac{a_1}{b}
\rfloor<\lfloor
\frac{a_2}{b}
\rfloor<
\cdots<\lfloor \frac{a_{k/4}}{b} \rfloor$ be the associated
$b$-floor wave.
Note that this $b$-floor wave is either an ascending wave or
an almost ascending wave with last difference  at most
$\frac{k/4}{10^{14}}$.  Hence, by Proposition 3.1, the number
of such
$b$-floor waves is at most, for $k$ sufficiently large,
$$
 k^2 \cdot 2^{\frac{13k}{100}} \cdot
\left(\frac{3}{2}\right)^{k/400}
\leq 2^{\frac{14k}{100}} \cdot
\left(\frac{3}{2}\right)^{k/400}
$$
(we have less than $k^2$ choices for $\lfloor \frac{a_1}{b}\rfloor$).

 Note that $A$ is monochromatic of
color, say $c$, only if none of the blocks $B_{\lfloor \frac{a_i}{b}
\rfloor}$, $1 \leq i \leq \frac{k}{4}$, is colored by $\gamma_c$.
Note that all of these blocks are at least $6 (=2r)$ blocks from
each other. We use Lemma 2.1 to give us that the probability that
$A$ is monochromatic is no more than
$$
3\left(\frac{2}{3}\right)^{k/4}.
$$

Thus, the probability that in a random $3$-coloring
of $[1,Q]$ we have a monochromatic ascending wave with
last difference less than $\frac{kb}{4 \cdot 10^{14}}$ is at most
$$
3\cdot 2^{\frac{14k}{100}} \cdot
\left(\frac{2}{3}\right)^{99k/400}.
$$
We have $3 < \left(\frac{3}{2}\right)^{3k/400}$ for $k$ sufficiently
large, so that the above probability is less than
$$
2^{\frac{14k}{100}} \cdot
\left(\frac{2}{3}\right)^{24k/100}.
$$
The above quantity is,
in particular, less than
$1/2$ for $k$ sufficiently large.
\B

To finish proving Theorem 3, we  apply
 Corollary 2.4 and Lemma 3.2 to
show that, for $k$ sufficiently large,
 there exists a $3$-coloring
of
$[1,Q]$ such that both of the following hold:

\ni
1) Any $\frac{k}{2}$-term monochromatic ascending wave has
last difference greater than $6b (=2rb)$.

\ni
2) Any $\frac{k}{4}$-term monochromatic ascending wave with
first difference greater than $6b (=2rb)$ has last difference
greater than $\frac{kb}{4 \cdot 10^{14}}$.

Hence, we conclude that there is a $3$-coloring of $[1,Q]$
such that any $\frac{3k}{4}$-term monochromatic ascending
wave has last difference greater
than $\frac{kb}{4 \cdot 10^{14}}$, for
$k$ sufficiently large.
This implies that the last term of such a monochromatic
ascending wave must be at least
$
\frac{3k}{4}+\frac{kb}{4 \cdot 10^{14}} \cdot
\frac{k}{4}.
$

We have $b=AW\left(\frac{k}{10(4r-4)};r-1\right)-1$
with $r=3$.  By 
Alon and Spencer's result [1], this gives us
$$
b \geq \frac{k^3}{10^{25} \cdot 8^{3}}
$$
for $k$ sufficiently large.

Hence, for $k$ sufficiently large, the last term
of an ascending wave
must be at least
$$
\frac{3k}{4}+\frac{k^2}{4^2 \cdot 10^{14}} \cdot \frac{k^3}{10^{25}
\cdot 8^{3}} >  \frac{k^5}{2^{13} \cdot 10^{39} } = Q.
$$

Since we have the existence of a $3$-coloring of $[1,Q]$
with no monochromatic $k$-term ascending wave, this
completes the proof of Theorem 3.

\ni
{\bf Remark}  From the lower bound
given in Proposition 3.1, it is not possible
to show
that $AW(k;r) \geq O(k^{2r-1})$ for
$r \geq 4$, which we believe is correct,
by using the argument presented in
Sections 2 and 3.

\section*{\normalsize References}
\footnotesize
\parindent=0pt

[1] N. Alon and J. Spencer, Ascending waves,
{\it J. Combin. Theory, Series A} {\bf 52} (1989),
275-287.

[2] T. Brown, P. Erd\H{o}s, and A. Freedman, Quasi-progressions
and descending waves, {\it J. Combin. Theory Series A}
{\bf 53} (1990), 81-95.

[3] J. Conway and R. Guy, The Book of Numbers, {\it Springer-Verlag},
New York, p. 106, 1996.

[4] B. Landman and A. Robertson, Ramsey Theory on the Integers,
{\it American Math. Society}, Providence, RI, 317pp., 2003.

[5] H. Robbins, A remark on Stirling's formula,
{\it American Math. Monthly} {\bf 62} (1955), 26-29.

\end{document}